\documentclass[a4paper,11pt]{article}
\usepackage[utf8x]{inputenc}
\usepackage{amsthm,amsmath,amssymb,oldgerm}
\usepackage[margin=3cm]{geometry}
\theoremstyle{plain}
\newtheorem{thm}{Theorem}

\newtheorem{cor}[thm]{Corollary}
\theoremstyle{definition}
 
\theoremstyle{definition}
\newtheorem{rem}[thm]{Remark}

\let\Im\relax
\DeclareMathOperator{\Im}{Im}
\DeclareMathOperator{\bk} {{\it{\mathcal{B}_{X}^{k}}}}
\DeclareMathOperator{\bkv} {{\it{\mathcal{B}_{X}^{k,\nu}}}}

\DeclareMathOperator{\bnkvn} {{\it{\mathcal{B}_{X}^{nk,\nu^n}}}}
\DeclareMathOperator{\bktwo} {{\it{\mathcal{B}_{X}^{2}}}}
\DeclareMathOperator{\bkm} {{\it{\mathcal{B}_{M,\mathcal{L}}^{k}}}}
\DeclareMathOperator{\hkm} {{\it{K_{M,\mathcal{L}}^{k}}}}

\DeclareMathOperator{\hyp}{\mu_{hyp}}

\DeclareMathOperator{\shyp}{\mu_{shyp}}
\let\Re\relax
\DeclareMathOperator{\Re}{Re}
\DeclareMathOperator{\vx}{\mathrm{vol_{\mathrm{hyp}}}}

\DeclareMathOperator{\lk}{{\it{\overline{\square}_{k}}}}
\let\ell\relax
\DeclareMathOperator{\ell}{\mathcal{L}}
\title{Heat kernels, Bergman kernels, and cusp forms}
{\small\author{Anilatmaja Aryasomayajula}}
\date{}
\begin{document}
\maketitle
\begin{abstract}
\noindent 
In this article, we describe a geometric method to study cusp forms, which relies on heat kernel and Bergman 
kernel analysis. This new approach of applying techniques coming from analytic geometry is based on the micro-local 
analysis of the heat kernel and the Bergman kernel from \cite{bouche} and \cite{berman}, respectively, 
using which we derive sup-norm bounds for cusp forms of integral weight, half-integral weight, and real weight 
associated to a Fuchsian subgroup of first kind. 

\vspace{0.2cm}\noindent
Mathematics Subject Classification (2010): 11F11, 11F37, 30C40.
\end{abstract}
\section{Introduction}\label{introduction}
This is both a survey and a research article elucidating heat kernel and Bergman kernel methods for studying cusp forms. We describe a 
geometric approach of Bouche and Berman to study sup-norm bounds for sections of a positive 
line bundle defined over a compact complex manifold. We then apply these methods to study cusp forms associated to 
a Fuchsian subgroup of first kind, which yields optimal results when the Fuchsian subgroup is cocompact. However, 
this approach  does not give optimal results when the Fuchsian subgroup is cofinite. But an extension of the 
methods of Bouche to cuspidal neighborhoods should allow one to derive optimal bounds for 
cusp forms, even when the Fuchsian subgroup is cofinite. 

\vspace{0.2cm}
\subsection{Notation}
Let $\mathbb{C}$ denote the complex plane. For $z\in\mathbb{C}$, let $x=\Re(z)$ and $y=\Im(z)$ denote the real and 
imaginary parts of $z$, respectively. Let 
\begin{equation*}
 \mathbb{H}=\lbrace  z\in \mathbb{C}|\,y= \Im(z) > 0 \rbrace
\end{equation*}
be the upper half-plane. Let $\Gamma \subset \mathrm{PSL}_{2}(\mathbb{R})$ be a Fuchsian subgroup of the 
first kind acting by fractional linear transformations on $\mathbb{H}$. We assume that $\Gamma$ admits no elliptic 
elements. 

\vspace{0.2cm}
Let $X$ be the quotient space $\Gamma\backslash \mathbb{H}$ of genus $g>1$. The quotient space $X$ admits the 
structure of a hyperbolic Riemann surface of finite volume. We allow $X$ to have genus $g=1$, if $X$ is not compact.  

\vspace{0.2cm}
We denote the $(1,1)$-form corresponding to the hyperbolic metric of $X$, which is compatible with the 
complex structure on $X$ and has constant negative curvature equal to minus one, by $\hyp(z)$. Locally, for 
$z\in X$, it is given by
\begin{equation*}
 \hyp(z)= \frac{i}{2}\cdot\frac{dz\wedge d\overline{z}}{{\Im(z)}^{2}}.
\end{equation*} 
Let $\shyp(z)$ denote the rescaled hyperbolic metric $\hyp(z)\slash \vx(X)$, which measures the volume of $X$ to be 
one. 

\vspace{0.2cm}
For $k\in\mathbb{R}_{>0}$, let $\nu$ denote the factor of automorphy of weight $k$ with the associated character 
being unitary. Let $S^{k}(\Gamma,\nu)$ denote the complex vector space of weight-$k$ cusp forms with respect to 
$\Gamma$ and $\nu$. Let $\lbrace f_{1},\ldots,f_{j_{k}} \rbrace$ denote an orthonormal basis of $S^{k}(\Gamma,\nu)$ 
with respect to the Petersson inner product. Then, for $z\in X$, put 
\begin{align*}
\bkv(z):=\sum_{i=1}^{j_{k}} y^{k}|f_{i}(z)|^{2}.
\end{align*}
When the associated character is trivial, we put $\bkv(z)=\bk(z)$.

\subsection{Sup norm bounds for the function $\bkv(z)$}
Let the Fuchsian subgroup $\Gamma$ be cocompact, i.e. $X$ is a compact Riemann surface. With notation as above, 
for $k\in\frac{1}{2}\mathbb{Z}$ (or $2\mathbb{Z}$), we have the following estimate
\begin{align}\label{estimate1}
\lim_{k}\sup_{z\in X}\frac{1}{k}\bk(z)= O(1),
\end{align}
where the implied constant is independent of $\Gamma$. 

\vspace{0.2cm}
Furthermore, for a fixed $k\in\mathbb{R}_{>0}$, let $\nu$ denote the factor of automorphy of weight $k$ with the 
associated character being unitary. Then, with notation as above, we have the following estimate
\begin{align}\label{estimate2}
\lim_{n}\sup_{z\in X}\frac{1}{nk}\bnkvn(z)= O(1),
\end{align}
where $n\in\mathbb{Z}$, and the implied constant is independent of $\Gamma$. 

\vspace{0.2cm}
Let $\Gamma$ now be a cofinite subgroup, i.e., $X$ is a noncompact hyperbolic Riemann surface of finite volume, and 
let $A$ be a compact subset of $X$. Then, with notation as above, for $k\in\frac{1}{2}\mathbb{Z}$ (or $2\mathbb{Z}$) 
and $z\in A$, we have the following estimate
\begin{align}\label{estimate3}
\lim_{k}\frac{1}{k}\bk(z)= O_{A}(1),
\end{align}
where the implied constant depends on $A$.

\vspace{0.2cm}
Furthermore, for a fixed $k\in\mathbb{R}_{>0}$, let $\nu$ denote the factor of automorphy of weight $k$ with the 
associated character being unitary. Let $A$ be any compact subset of $X$. Then, with notation as above, 
for any $z\in A$, we have the following estimate
\begin{align}\label{estimate4}
\lim_{n}\frac{1}{nk}\bnkvn(z)= O_{A}(1),
\end{align}
where $n\in\mathbb{Z}$ and the implied constant depends on $A$.

\vspace{0.2cm}
Our estimates \eqref{estimate1}, \eqref{estimate2}, \eqref{estimate3}, and \eqref{estimate4} are optimal. 
However, when $X$ is noncompact we cannot extend estimates \eqref{estimate3} and \eqref{estimate4} to the entire 
Riemann surface, i.e., our method does not yield optimal estimates when $A$ is equal to $X$. However, an extension of 
Bouche's methods to cuspidal neighborhoods will enable the extension of estimates \eqref{estimate3} and 
\eqref{estimate4} to $X$.

\vspace{0.2cm}
Lastly, our methods extend with notational changes to higher dimensions, namely to Hilbert modular cusp forms and 
Siegel modular cusp forms. 

\subsection{Existing results on sup-norm bounds for the function $\bkv(z)$}

In \cite{jk2}, using heat kernel analysis, Jorgenson and Kramer derived sup-norm bounds for the the Bergman kernel $\bktwo(z)$, 
associated to any hyperbolic Riemann surface $X$ (compact or noncompact of finite volume). The bounds of Jorgenson and Kramer are 
optimal. Especially for the case $X=Y_{0}(N)$, they derived
\begin{align*}
\sup_{z\in Y_{0}(N)}\mathcal{B}_{Y_{0}(N)}^{2}(z)=O(1),
\end{align*}
where the implied constant does not depend on the modular curve $Y_{0}(N)$. 

\vspace{0.2cm}
In \cite{jk1}, extending their method from \cite{jk2}, Jorgenson, Kramer, and Friedman derived sup-norm bounds for the 
the Bergman kernel $\bk(z)$, associated to any hyperbolic Riemann surface $X$ (compact or noncompact of finite 
volume). When $X$ is a compact hyperbolic Riemann surface, they showed that
\begin{align*}
 \sup_{z\in X}\bk(z)=O(k),
\end{align*}
where the implied constant is independent of the Riemann surface $X$. When $X$ is a noncompact hyperbolic Riemann 
surface of finite volume, they showed that 
\begin{align*}
 \sup_{z\in X}\bk(z)=O(k^{\frac{3}{2}}),
\end{align*}
where the implied constant is independent of the Riemann surface $X$. The estimates of Jorgenson and Kramer are optimal, 
as shown in \cite{jk2}.

\vspace{0.2cm}
It is possible to extend the heat kernel analysis of Jorgenson and Kramer to higher dimensions, namely to Hilbert modular 
cusp forms and Siegel modular cusp forms of both integral and half-integral weight. However, one has to address 
certain non trivial convergence issues, while doing so. 

\vspace{0.2cm}
For $k\in\frac{1}{2}\mathbb{Z}$ and $N\in\mathbb{N}$, let $f$ be any weight-$k$ cusp form with respect to the arithmetic subgroup 
$\Gamma_{0}(4N)$. Furthermore, let $f$ be normalized with respect to the Petersson 
inner-product. Then, in \cite{kiral}, Kiral has derived the following estimate 
\begin{align*}
 \sup_{z\in Y_{0}(N)}y^{k}|f(z)|^{2}=O_{k,\varepsilon}\big(N^{\frac{1}{2}-\frac{1}{18}+\varepsilon}\big),
\end{align*}
for any $\varepsilon>0$. Using above estimate, one can derive 
\begin{align*}
 \sup_{z\in Y_{0}(N)}\mathcal{B}_{Y_{0}(N)}(z)=O_{k,\varepsilon}\big(N^{\frac{1}{2}-\frac{1}{18}+\varepsilon}\big),
\end{align*}
for any $\varepsilon>0$. 

\vspace{0.2cm} 
For $k\in\mathbb{R}_{>0}$ with $k>2$, the Bergman kernel $\bk(z)$ can be represented by an infinite series, which is 
uniformly convergent in $z\in X$. Using which, Steiner has extended the bounds of Jorgenson and Kramer to real weights. 

\vspace{0.2cm} 
Let $\Gamma$ any subgroup of finite index in $\mathrm{SL}_{2}(\mathbb{R})$, and for a fixed $k\in\mathbb{R}_{>0}$ with 
$k\gg1$, let $\nu$ denote the factor of automorphy of weight $k$ with the associated character being unitary. 
Furthermore, let $A$ be a compact subset of $X$. Then, in \cite{raphael}, Steiner has derived the following estimates
\begin{align*}
 \sup_{z\in A}\bkv(z)=O_{A}(k),
\end{align*}
where the implied constant depends on the compact subset $A$; and 
\begin{align*}
 \sup_{z\in X }\bkv(z)=O_{X}\big(k^{\frac{3}{2}}\big),
\end{align*}
where the implied constant depends on $X$. 

\section{Heat kernels and Bergman kernels on compact complex manifolds}
In this section, we recall the main results from \cite{bouche} and \cite{berman}, which we use in 
the next section. 

\vspace{0.2cm}
Let $(M,\omega)$ be a compact complex manifold of dimension $n$ with a Hermitian metric $\omega$. Let 
$\ell$ be a positive Hermitian holomorphic line bundle on $M$ with the Hermitian metric given by 
$\|s(z)\|^{2}_{\ell}:=e^{-\phi(z)}|s(z)|^{2}$, where $s\in\ell$ is any section, and $\phi(z)$ is 
a real-valued function defined on $M$. 

\vspace{0.2cm}
For any $k\in\mathbb{N}$, let $\lk:=(\overline{\partial}^{\ast}+\overline{\partial})^{2}$ denote the $\overline{\partial}$-Laplacian 
acting on smooth sections of the line bundle $\ell^{\otimes k}$. Let $\hkm(t;z,w)$ denote the smooth kernel of the 
operator $e^{-\frac{2t}{k}\lk}$. We refer the reader to p. 2 in \cite{bouche}, for the details regarding the 
properties which uniquely characterize the heat kernel $\hkm(t;z,w)$. When $z=w$, the heat kernel $\hkm(t;z,w)$ 
admits the following spectral expansion
\begin{align}\label{spectralexpn}
\hkm(t;z,w)=\sum_{n\geq 0} e^{-\frac{2t}{k}\lambda_{n}^{k}}\varphi_{n}(z)\otimes\varphi_{n}^{\ast}(w),
\end{align}
where $\lbrace\lambda_{n}^k\rbrace_{n\geq0}$ denotes the set of eigenvalues of $\lk$ 
(counted with multiplicities), and $\lbrace\varphi_{n}\rbrace_{n\geq0}$ denotes a set of associated 
orthonormal eigenfunctions. 

\vspace{0.2cm}
Let $\lbrace s_{i}\rbrace$ denote an orthonormal basis of $H^{0}(M,\mathcal{L}^{\otimes k})$. For any $z\in M$, 
the Bergman kernel is given by
\begin{align}\label{bkdefn}
\bkm(z):= \sum_{i}\| s_{i}(z)\|_{\mathcal{L}^{\otimes k}}^{2}.
\end{align}
For any $z\in M$ and $t\in\mathbb{R}_{>0}$, from the spectral expansion of the heat kernel $\hkm(t;z,w)$ described in equation \eqref{spectralexpn}, 
it is easy to see that
\begin{align}\label{hkbkreln}
\bkm(t;z)\leq \hkm(t;z,z)\quad \mathrm{and}\quad
\lim_{t} \hkm(t;z,z)=\bkm(t;z).
\end{align}
For $z\in M$, let $c_{1}(\ell)(z):=\frac{i}{2\pi}\partial\overline{\partial}\phi(z)$ denote the first Chern form 
of the line bundle $\ell$. Let $\alpha_{1},\ldots,\alpha_{n}$ denote the eigenvalues of $\partial\overline{\partial}\phi(z)$ at 
the point $z\in M$. Then, with notation as above, 
from Theorem 1.1 in \cite{bouche}, for any $z\in M$ and $t\in (0,k^{\varepsilon})$, and for 
a given $\varepsilon>0$ not depending on $k$, we have
\begin{align}\label{boucheeqn1}
\lim_{k}\frac{1}{k^{n}}\hkm(t;z,z)=\prod_{j=1}^{n}\frac{\alpha_{j}}{(4\pi)^{n}\sinh(\alpha_{j}t)}, 
\end{align}
and the convergence of the above limit is uniform in $z$. 

\vspace{0.2cm}
Using equations \eqref{hkbkreln} and \eqref{boucheeqn1}, in Theorem 2.1 in \cite{bouche}, Bouche derived the 
following asymptotic estimate 
\begin{align}\label{boucheeqn2}
\lim_{k}\frac{1}{k^{n}}\bkm(z)= O\big(\mathrm{det}_{\omega}\big(c_{1}(\ell)(z)\big)\big),
\end{align}
where the implied constant does not depend on $X$, and the convergence of the above limit is uniform in $z\in X$. 

\vspace{0.2cm}
When $M$ is a noncompact complex manifold, using micro-local analysis of the Bergman kernel, in \cite{berman}, 
Berman derived the following estimate
\begin{align*}
\limsup_{k}\frac{1}{k^{n}}\bkm(z)\leq \mathrm{det}_{\omega}\big(c_{1}(\ell)(z)\big).
\end{align*}

\vspace{0.2cm}
Furthermore, let $A$  be any compact subset of $M$. Then, for any $z\in A$, from the proof of Corollary 3.3 in 
\cite{berman}, we have
\begin{align}\label{bermaneqn}
\lim_{k}\frac{1}{k^{n}}\bkm(z)=O_{A}\bigg(\mathrm{det}_{\omega}\big(c_{1}(\ell)(z)\big)\bigg),
\end{align}
where the implied constant depends on the compact subset $A$.  
\section{Estimates of cusp forms}
In this section, using results from previous section, we prove estimates \eqref{estimate1}, \eqref{estimate2}, 
\eqref{estimate3}, and \eqref{estimate4}.

\vspace{0.2cm}
Let notation be as in Section \ref{introduction}. Let $\Omega_{X}$ denote the cotangent bundle over $X$. 
Then, for any $k\in2\mathbb{Z}$, cusp forms of 
weight $k$ with respect to $\Gamma$ are global section of the line bundle $\Omega_{X}^{\otimes k\slash 2}$. 
Furthermore, recall that for any $f\in\Omega_{X}$, i.e., $f$ a weight-$2$ cusp form, the Petersson metric on the 
line bundle $\Omega_{X}$ is given by 
\begin{align}\label{peterssonip1}
 \|f(z)\|_{\Omega_{X}}^{2}:=y^{2}|f(z)|^{2}. 
\end{align}
Let $\omega_{X}$ denote the line bundle of cusp forms of weight $\frac{1}{2}$ over $X$. Then, for any 
$k\in\frac{1}{2}\mathbb{Z}$, cusp forms of weight-$k$ with respect to $\Gamma$ are global section of the line 
bundle $\omega_{X}^{\otimes 2k}$. Furthermore, recall that for any 
$f\in\omega_{X}$, i.e., $f$ a weight-$\frac{1}{2}$ cusp form, the Petersson metric on the line 
bundle $\omega_{X}$ is given by 
\begin{align}\label{peterssonip2}
 \|f(z)\|_{\omega_{X}}^{2}:=y^{1\slash 2}|f(z)|^{2}. 
\end{align}

\begin{rem}
For any $z\in X$ and $k\in 2\mathbb{Z}$, from the definition of the Bergman kernel $\mathcal{B}_{X,\Omega_{X}}^{k
\slash 2}(z)$ for the line bundle $\Omega_{X}^{\otimes k\slash 2}$ from equation \eqref{bkdefn}, we have
\begin{align*}
\mathcal{B}_{X,\Omega_{X}}^{k\slash 2}(z)=\bk(z). 
\end{align*}
Similarly, for any $z\in X$ and $k\in \frac{1}{2}\mathbb{Z}$, from the definition of the Bergman kernel 
$\mathcal{B}_{X,\Omega_{X}}^{2k}(z)$ for the line bundle $\omega_{X}^{\otimes 2k}$ from equation \eqref{bkdefn}, we have
\begin{align}\label{remeqn2}
\mathcal{B}_{X,\omega_{X}}^{2k}(z)=\bk(z). 
\end{align}
\end{rem}

\vspace{0.2cm}
\begin{thm}\label{thm1}
Let $\Gamma$ be cocompact, i.e., $X$ is a compact hyperbolic Riemann surface. Then, with notation as above, 
for $k\in\frac{1}{2}\mathbb{Z}$ (or $2\mathbb{Z}$), we have the following estimate
\begin{align*}
\lim_{k}\sup_{z\in X}\frac{1}{k}\bk(z)=O(1), 
\end{align*}
where the implied constant is independent of $\Gamma$.   
\begin{proof}
We refer the reader to Theorem 2 in \cite{anil} for the proof of the theorem. We briefly describe the proof of 
the theorem for $k\in \frac{1}{2}\mathbb{Z}$, and the case for $k\in2\mathbb{Z}$, follows automatically with 
notational changes. For any $z\in X$,  observe that
\begin{align*}
c_{1}(\omega_{X}^{\otimes 2})(z)=\frac{1}{4\pi}\hyp(z),
\end{align*}
which shows that the line bundle $\omega_{X}^{\otimes 2}$ is positive, and $\mathrm{det}_{\hyp}\big(c_{1}(
\omega_X^{\otimes 2})(z)\big)=\frac{1}{4\pi}$. Using equation \eqref{remeqn2}, 
and applying estimate \eqref{boucheeqn2} to the complex manifold $X$ with its natural 
Hermitian metric $\hyp$ and the line bundle $\omega_{k}^{\otimes 2k}$, we find
\begin{align*}
\lim_{k}\frac{1}{k}\bk(z)=\lim_{k}\frac{1}{k}\mathcal{B}_{X,\omega_X}^{2k}(z)=O
\bigg(\mathrm{det}_{\hyp}\big(c_{1}(\omega_X^{\otimes 2})(z)\big)\bigg)=O(1).
\end{align*}
As the above limit convergences uniformly in $z\in X$, and as $X$ is compact, we have
\begin{align*}
\sup_{z\in X} \lim_{k}\frac{1}{k}\bk(z)=\lim_{k}\sup_{z\in X}\frac{1}{k}\bk(z)=O(1),
\end{align*}
which completes the proof of the theorem.
\end{proof}
\end{thm}

\begin{cor}
Let $\Gamma$ be cofinite, i.e., $X$ is a noncompact hyperbolic Riemann surface of finite volume, and let $A$ be a 
compact subset of $X$. Then, with notation as above, for $k\in\frac{1}{2}\mathbb{Z}$ (or $2\mathbb{Z}$), and 
$z\in A$, we have the following estimate
\begin{align*}
\lim_{k}\frac{1}{k}\bk(z)= O_{A}(1),
\end{align*}
where the implied constant depends on $A$.
\begin{proof}
The proof of the theorem follows from estimate \eqref{bermaneqn}, and from similar arguments as in Theorem \ref{thm1}. 
 \end{proof}
\end{cor}

\begin{rem}
For a fixed $k\in\mathbb{R}_{>0}$, and let $\omega_{X,k,\nu}$ denote the line bundle of weight-$k$ cusp forms with 
the factor of automorphy $\nu$, and associated character being unitary. Then, for any $n\in\mathbb{Z}$, cusp forms of 
weight-$nk$ with respect to $\Gamma$ and $\nu^{n}$ are global section of the line 
bundle $\omega_{X,k,\nu^n}^{\otimes n}$. Furthermore, recall that for any 
$f\in\omega_{X,k,\nu}$, the Petersson metric on the line bundle $\omega_{X,k,\nu}$ is given by 
\begin{align}\label{peterssonip3}
 \|f(z)\|_{\omega_{X,k,\nu}}^2:=y^{k}|f(z)|^{2}. 
\end{align}
For any $z\in X$ and $n\in \mathbb{Z}$, from the definition of the Bergman kernel $\mathcal{B}_{X,\omega_{X,k,\nu^n}}^{n}
(z)$ for the line bundle $\omega_{X,k,\nu^n}^{\otimes n}$ from equation \eqref{bkdefn}, we have
\begin{align}\label{remeqn3}
\mathcal{B}_{X,\omega_{X,k,\nu^n}}^{n}(z)=\bnkvn(z). 
\end{align}
\end{rem}

\begin{thm}\label{thm2}
Let $\Gamma$ be cocompact, i.e., $X$ is a compact hyperbolic Riemann surface. For a fixed $k\in\mathbb{R}_{>0}$, 
let $\nu$ denote the factor of automorphy of weight $k$ with the associated character being unitary. Then, with 
notation as above, we have the following estimate
\begin{align*}
\lim_{n}\sup_{z\in X}\frac{1}{nk}\bnkvn(z)= O(1),
\end{align*}
where $n\in\mathbb{Z}$, and the implied constant is independent of $\Gamma$. 
\begin{proof}
From equation \eqref{peterssonip3}, for any $z\in X$, observe that
\begin{align*}
c_{1}\big(\omega_{X,k,\nu}\big)(z)=-\frac{i}{2\pi}\partial\overline{\partial}\log\big(y^{k}|f(z)|^{2}\big)
=\frac{k}{4\pi}\hyp(z),
\end{align*}
which shows that the line bundle $\omega_{X,k,\nu}$ is positive, and $\mathrm{det}_{\hyp}\big(c_{1}
(\omega_{X,k,\nu})(z)\big)=\frac{k}{4\pi}$. Using equation \eqref{remeqn3}, and applying estimate 
\eqref{boucheeqn2} to the complex manifold $X$ with its natural Hermitian metric $\hyp$ and the line bundle 
$\omega_{X,k,\nu^n}^{\otimes n}$, we find
\begin{align*}
\lim_{n}\frac{1}{nk}\bnkvn(z)=\lim_{n}\frac{1}{nk}\mathcal{B}_{X,\omega_{X,k,\nu^n}}^{n}(z)=O
\bigg(\frac{1}{k}\mathrm{det}_{\hyp}\big(c_{1}\big(\omega_{X,k,\nu}\big)(z)\big)\bigg)=O(1).
\end{align*}
As the above limit convergences uniformly in $z\in X$, and as $X$ is compact, we have
\begin{align*}
\sup_{z\in X} \lim_{n}\frac{1}{nk}\bnkvn(z)=\lim_{n}\sup_{z\in X}\frac{1}{nk}\bnkvn(z)=O(1),
\end{align*}
which completes the proof of the theorem.
\end{proof}
\end{thm}

\begin{cor}\label{cor1}
Let $\Gamma$ be cofinite, i.e., $X$ is a noncompact hyperbolic Riemann surface of finite volume. For a fixed 
$k\in\mathbb{R}_{>0}$, let $\nu$ denote the factor of automorphy of weight $k$ with the associated character 
being unitary. Let $A$ be a compact subset of $X$. Then, with notation as above, 
for any $z\in A$, we have the following estimate
\begin{align*}
\lim_{nk}\frac{1}{nk}\bnkvn(z)= O_{A}(1),
\end{align*}
where $n\in\mathbb{Z}$, and the implied constant depends on $A$.
\begin{proof}
The proof of the theorem follows from estimate \eqref{bermaneqn}, and from similar arguments as in Theorem \ref{thm2}. 
\end{proof}
\end{cor}
\paragraph{Acknowledgements}
The author would like to thank J.~Kramer and J.~Jorgenson for introducing him to the area of automorphic forms 
and heat kernels. The author would like to express his gratitude to T. Bouche for providing him with helpful 
references, and to Archana S. Morye, for many helpful discussions and remarks. The author would also like to thank 
R. S. Steiner for his communications.  

\vspace{0.3cm}
Anilatmaja Aryasomayajula\\ {\it{anilatmaja@gmail.com}}  \\ 
Department of Mathematics, \hfill\\University of Hyderabad, \\Prof. C.~R.~Rao Road, Gachibowli,\\
Hyderabad, 500046, India.

\begin{thebibliography}{99}
\addcontentsline{toc}{chapter}{Bibliography}
\bibitem{anil} A. Aryasomayajula, Heat kernel approach for sup-norm bounds for cusp forms of integral and half-integral 
weight, arxiv preprint, arXiv:1506.08497.
\bibitem{berman} R. J. Berman, Bergman kernels and local holomorphic Morse inequalities, Math. Z. 248, 
(2004), 325--344.
\bibitem{bouche} T.~Bouche, Asymptotic results for Hermitian line bundles over complex manifolds: The heat kernel 
approach, Higher-dimensional complex varieties, 67--81, de Gruyter, Berlin, 1996.
\bibitem{kiral} E. M. Kiral, Bounds on Sup-norms of Half Integral Weight Modular Forms, Acta Arithmetica, Vol. 165, (2014), 385-399.
\bibitem{jk1} J.~Friedman, J.~Jorgenson, and J.~Kramer, Uniform sup-norm bounds on average for
cusp forms of higher weights, arxiv preprint, arXiv:1305.1348 .
\bibitem{jk2} J.~Jorgenson and J.~Kramer, Bounding the sup-norm of automorphic
forms, GAFA 14, (2004), 1267--1277.
\bibitem{raphael} R. S. Steiner, Uniform bounds on sup-norms of holomorphic forms of real weight, arxiv preprint, 
arXiv:1406.2918.
\end{thebibliography}
\end{document}